\documentclass{article}
\usepackage{amssymb}
\usepackage{amsfonts}

\input{tcilatex}
\begin{document}

\begin{center}
\bigskip \textbf{The boundary conrollability for abstract wave equations and
its applications}

\textbf{Veli\ B. Shakhmurov}

Department of Mechanical Engineering, Istanbul Okan University, Akfirat,
Tuzla 34959 Istanbul, E-mail: veli.sahmurov@okan.edu.tr,

Institute of Mathematics and Mechanics, Azerbaijan National Academy of
Sciences, AZ1141, Baku, F. Agaev 9

E-mail: veli.sahmurov@gmail.com

\bigskip \textbf{Abstract}

The paper is devoted to the exact controllability of a system of coupled
abstract wave equations when the control is exerted on a part of the
boundary by means of one control. We give a Kalman \ type condition and give
a description of the attainable set.

\bigskip \textbf{1. Introduction, definitions}
\end{center}

We consider here, the controllability properties of the nonlocal mixed
problem for abstract wave equation%
\begin{equation}
u_{tt}-u_{xx}+Au=0\text{, }\left( x,t\right) \in \text{ }Q=\left( 0,a\right)
\times \left( 0,T\right) ,  \tag{1.1}
\end{equation}%
\[
\alpha _{1}u\left( 0,t\right) +\beta _{1}u\left( a,t\right) =bf\left(
t\right) ,\text{ }\alpha _{2}u\left( 0,t\right) +\beta _{2}u\left(
a,t\right) =0\text{ for }t\in \left( 0,T\right) , 
\]%
\[
u\left( x,0\right) =u_{0}\left( x\right) ,\text{ }u_{t}\left( x,0\right)
=u_{1}\left( x\right) \text{ for }x\in \left( 0,a\right) , 
\]%
where $a,T>0$ are given numbers, $\alpha _{i},\beta _{i}$ are given
generally complex numbers, $A$ is a linear operator in a Hilbert space $H,$ $%
b$ is a given element in $H$ and $f\in L^{2}\left( 0,T\right) $ is a control
function to be determined which acts on the equation by means of the
nonlocal boundary condition $\left( 1.1\right) $. The initial data $\left(
u_{0},u_{1}\right) $ will belong to a Hilbert space $\mathbb{H=H}_{0}\times 
\mathbb{H}_{1}$, where $\mathbb{H}_{0}\subset \mathbb{H}_{1}$ which is to be
specified in our main result. Our goal is to give suficient conditions for
the exact boundary controllability of the problem $\left( 1.1\right) $ by
using the given space $\mathbb{H}$.

We recall that the problem $\left( 1.1\right) $ is exactly controllable in $%
\mathbb{H}$ at time $T$ if, for every initial and final data $\left(
u_{0},u_{1}\right) $, $\left( \upsilon _{0},\upsilon _{1}\right) $ both in $%
\mathbb{H}$, there exists a control $f\in L^{2}\left( 0,T\right) $ such that
the solution of problem $\left( 1.1\right) $ corresponding to $\left(
u_{0},u_{1},f\right) $ satisfies 
\begin{equation}
u\left( x,T\right) =\upsilon _{0}\left( x\right) ,\text{ }u_{t}\left(
x,T\right) =\upsilon _{1}\left( x\right) \text{ for }x\in \left( 0,a\right) .
\tag{1.2}
\end{equation}

Due to the linearity and time reversibility of problem $\left( 1.1\right) $,
this is equivalent to exact controllability from zero at time $T$.

The controllability properties of problem $\left( 1.1\right) $\ are well
known for $d-$dimensional coupled wave equations, i.e. for case of $H=%
\mathbb{R}^{d}$, $A$ is is a given coupling matrix and $b$ is a given vector
from $\mathbb{R}^{d}$ (see e.g. \ $\left[ 1-12\right] $). Most of the known
controllability results of $\left( 1.1\right) $ when $A$ is matix are in the
case of two coupled equations (see $\left[ \text{1, 12}\right] $). But the
results are for a particular coupling matrix $A$. In the $d$-dimensional
situation, that is, for a system of coupled wave equations in a domain $%
\mathbb{R}^{n}$, Alabau-Boussouria and collaborators have obtained several
results in the case of two equations with the Laplacian plus additional zero
order terms and particular coupling matrices (see e.g. $\left[ 1-3\right] $
and the references therein). On the other hand, controllability properties
of linear ordinary differential systems are well understood. In contrast to
these above mentioned workes, the problem $\left( 1.1\right) $ involves
generally infinite dimensionel linear operator $A$ in abstract Hilbert space 
$H$. Moreover, the boundary value problem (BVP) is nonlocal, in general. If
we put $\beta _{1}=\alpha _{2}=0$, then the nonlocal mixed problem $\left(
1.1\right) $ stated to be \ a local mixed problem. Under some sufficient
condition on $\alpha _{i}$, $\beta _{i}$, $i=1,2$, $b,f$ \ and operator $A$
we derive the conrollability of $\left( 1.1\right) .$

By selecting the space $H$ and the operators $A$ in $\left( 1.1\right) $, we
obtain different boundary controllability proplem with nonlocal conditions
for wave equations which occur in application. Let we put $H=l_{2}$ and
choose $A$\ as infinite matrices $\left[ a_{mj}\right] $ for $m,j=1,2,...,N,$
$N\in \mathbb{N},$ where $\mathbb{N-}$denote the set of natural numbers.
Then from our results we obtain the exact boundary controllability of the
mixed problem for infinite many system of wave equations 
\begin{equation}
\partial _{t}^{2}u_{m}-\partial
_{x}^{2}u_{m}+\dsum\limits_{j=1}^{m}a_{mj}u_{m}=0,\text{ }\left( x,t\right)
\in \text{ }Q,  \tag{1.3}
\end{equation}%
\[
\alpha _{1}u_{m}\left( 0,t\right) +\beta _{1}u_{m}\left( a,t\right)
=bf\left( t\right) ,\text{ }\alpha _{2}u_{m}\left( 0,t\right) +\beta
_{2}u_{m}\left( a,t\right) =0\text{ for }t\in \left( 0,T\right) , 
\]%
\[
u_{m}\left( x,0\right) =u_{m0}\left( x\right) ,\text{ }\partial
_{t}u_{m}\left( x,0\right) =u_{m1}\left( x\right) \text{ for }x\in \left(
0,a\right) , 
\]%
where $a_{mj}=a_{mj}\left( x\right) $ are complex valued functions and $%
u_{j}=u_{j}\left( x,t\right) .$

Moreover, let we choose $E=L^{2}\left( 0,1\right) $ and $A$ to be
degenerated differential operator in $L^{2}\left( 0,1\right) $ defined by 
\[
D\left( A\right) =\left\{ u\in W_{\gamma }^{\left[ 2\right] ,2}\left(
0,1\right) ,\right. \left. \nu _{k}u^{\left[ m_{k}\right] }\left( 0\right)
+\delta _{k}u^{\left[ m_{k}\right] }\left( 1\right) =0,\text{ }k=1,2\right\}
,\text{ } 
\]%
\begin{equation}
\text{ }A\left( x\right) u=a_{1}\left( x,y\right) u^{\left[ 2\right]
}+a_{2}\left( x,y\right) u^{\left[ 1\right] }\text{, }x\in \left( 0,a\right)
,\text{ }y\in \left( 0,1\right) ,\text{ }m_{k}\in \left\{ 0,1\right\} , 
\tag{1.4}
\end{equation}%
\ \ \ where $u^{\left[ i\right] }=\left( y^{\gamma }\frac{d}{dy}\right)
^{\gamma }u$ for $0\leq \gamma <\frac{1}{2},$ $a_{1}=b_{1}\left( x,y\right) $
is a cont\i nous, $a_{2}=b_{2}\left( x,y\right) $ is a bounded functon on $%
y\in $ $\left[ 0,1\right] $ for a.e. $x\in \left( 0,a\right) ,$ $\nu _{k}$, $%
\delta _{k}$ are complex numbers and $W_{\gamma }^{\left[ 2\right] ,2}\left(
0,1\right) $ is a weighted Sobolev spase defined by 
\[
W_{\gamma }^{\left[ 2\right] }\left( 0,1\right) =\left\{ {}\right. u:u\in
L^{2}\left( 0,1\right) ,\text{ }u^{\left[ 2\right] }\in L^{2}\left(
0,1\right) ,\text{ } 
\]%
\[
\left\Vert u\right\Vert _{W_{\gamma }^{\left[ 2\right] }}=\left\Vert
u\right\Vert _{L^{2}}+\left\Vert u^{\left[ 2\right] }\right\Vert
_{L^{2}}<\infty . 
\]%
Then, from $\left( 1.1\right) -\left( 1.2\right) $ we get\ the exact
boundary controllability of the mixed problem for idegenerate wave equations

\begin{equation}
u_{tt}-\partial _{x}^{2}u+\left( a_{1}\frac{\partial ^{\left[ 2\right] }u}{%
\partial y^{2}}+a_{2}\frac{\partial ^{\left[ 1\right] }u}{\partial y}\right)
=0,\text{ }  \tag{1.5}
\end{equation}%
\[
x\in \left( 0,a\right) ,\text{ }y\in \left( 0,1\right) ,\text{ }t\in \left(
0,T\right) ,\text{ }u=u\left( x,y,t\right) , 
\]%
\ \ \ 

\begin{equation}
\nu _{k}u^{\left[ m_{k}\right] }\left( x,0,t\right) +\delta _{k}u^{\left[
m_{k}\right] }\left( x,1,t\right) =0,\text{ }k=1,2,  \tag{1.6}
\end{equation}

\begin{equation}
u\left( x,y,0\right) =\varphi \left( x,y\right) ,\text{ }u_{t}\left(
x,y,0\right) =\psi \left( x,y\right) \text{.}  \tag{1.7}
\end{equation}

To state our results, we provide the following definition:

Let $E$ be a Banach space. $L^{p}\left( \Omega ;E\right) $ denotes the space
of strongly measurable $E$-valued functions that are defined on the
measurable subset $\Omega \subset \mathbb{R}^{n}$ with the norm

\[
\left\Vert f\right\Vert _{L^{p}}=\left\Vert f\right\Vert _{L^{p}\left(
\Omega ;E\right) }=\left( \int\limits_{\Omega }\left\Vert f\left( x\right)
\right\Vert _{E}^{p}dx\right) ^{\frac{1}{p}},\text{ }1\leq p<\infty \ . 
\]

Let $H$ be a Hilbert space and 
\[
\left\Vert u\right\Vert =\left\Vert u\right\Vert _{H}=\left( u,u\right)
_{H}^{\frac{1}{2}}\text{ for }u\in H. 
\]

For $p=2$ and $E=H$, $L^{p}\left( \Omega ;E\right) $ becomes a $H$-valued
function space with inner product: 
\[
\left( f,g\right) _{L^{2}\left( \Omega ;H\right) }=\int\limits_{\Omega
}\left( f\left( x\right) ,g\left( x\right) \right) _{H}dx,\text{ }f,\text{ }%
g\in L^{2}\left( \Omega ;H\right) . 
\]

Here, $W^{s,2}\left( \mathbb{R}^{n};H\right) $, $-\infty <s<\infty $ denotes
the $H-$valued Sobolev space of order $s$ which is defined as: 
\[
W^{s,2}=W^{s,2}\left( \mathbb{R}^{n};H\right) =\left( I-\Delta \right) ^{-%
\frac{s}{2}}L^{2}\left( \mathbb{R}^{n};H\right) 
\]%
with the norm 
\[
\left\Vert u\right\Vert _{W^{s,2}}=\left\Vert \left( I-\Delta \right) ^{%
\frac{s}{2}}u\right\Vert _{L^{2}\left( R^{n};H\right) }<\infty . 
\]%
It clear that $W^{0,2}\left( \mathbb{R}^{n};E\right) =L^{2}\left( \mathbb{R}%
^{n};H\right) .$ Let $H_{0}$ and $H$ be two Hilbert spaces and $H_{0}$ is
continuously and densely embedded into $H$. Let $W^{s,2}\left( \mathbb{R}%
^{n};H_{0},H\right) $ denote the Sobolev-Lions type space, i.e., 
\[
W^{s,2}\left( \mathbb{R}^{n};H_{0},H\right) =\left\{ u\in W^{s,2}\left( 
\mathbb{R}^{n};H\right) \cap L^{2}\left( \mathbb{R}^{n};H_{0}\right)
,\right. \text{ } 
\]%
\[
\left. \left\Vert u\right\Vert _{W^{s,2}\left( \mathbb{R}^{n};H_{0},H\right)
}=\left\Vert u\right\Vert _{L^{2}\left( \mathbb{R}^{n};H_{0}\right)
}+\left\Vert u\right\Vert _{W^{s,2}\left( \mathbb{R}^{n};H\right) }<\infty
\right\} . 
\]

Let $C\left( \Omega ;E\right) $ denote the space of $E-$valued uniformly
bounded continious functions on $\Omega $ with norm 
\[
\left\Vert u\right\Vert _{C\left( \Omega ;E\right) }=\sup\limits_{x\in
\Omega }\left\Vert u\left( x\right) \right\Vert _{E}. 
\]

$C^{m}\left( \Omega ;E\right) $\ will denote the spaces of $E$-valued
uniformly bounded strongly continuous and $m$-times continuously
differentiable functions on $\Omega $ with norm 
\[
\left\Vert u\right\Vert _{C^{m}\left( \Omega ;E\right) }=\max\limits_{0\leq
\left\vert \alpha \right\vert \leq m}\sup\limits_{x\in \Omega }\left\Vert
D^{\alpha }u\left( x\right) \right\Vert _{E}. 
\]

\textbf{Definition 1.1. }Let $S$ be a positive operator in a Hilbert space $%
H $ with discrete specrum $\left\{ \lambda _{k}\right\} _{k=1}^{\infty }$
and corresponding eigenfunctions $\left\{ \varphi _{k}\right\}
_{k=1}^{\infty }$. Moreover, we assume that $\left\{ \varphi _{k}\right\} $
is a ortonormal system in $H$. Let 
\[
l_{r}^{2}=\left\{ \left\{ c_{k}\right\} \text{: }\left\Vert \left\{
c_{k}\right\} \right\Vert _{l_{r}^{2}}=\left( \dsum\limits_{k=1}^{\infty
}\left\vert c_{k}\right\vert ^{2}\left\vert \lambda _{k}\right\vert
^{r}\right) ^{\frac{1}{r}}\right\} . 
\]

We then define the space 
\[
W_{r}=\left\{ f:f=\dsum\limits_{k=1}^{\infty }c_{k}\varphi _{k}\text{, }%
\left\Vert f\right\Vert =\left\Vert \left\{ c_{k}\right\} \right\Vert
_{l_{r}^{2}}<\infty \right\} . 
\]

For $r>0$, we set $W_{r}=D\left( S^{\frac{r}{2}}\right) $, where $D\left(
S^{\theta }\right) $ denotes the domain of the operator $S^{\theta }$. In
the case where $r=0$, $W_{r}$ $=H$, and for $r<0$, we set $W_{r}=\left(
W_{-r}\right) ^{\ast }$ where "$\ast $" indicates the dual space. Also, we
recall that the operator $-\partial _{x}^{2}$ with nonlocal boundary
conditions 
\begin{equation}
\alpha _{1}u\left( 0\right) +\beta _{1}u\left( a\right) =0,\text{ }\alpha
_{2}u\left( 0\right) +\beta _{2}u\left( a\right) =0\text{ }  \tag{1.8}
\end{equation}%
with $\eta =\alpha _{1}\beta _{2}-\alpha _{2}\beta _{1}\neq 0$\ admits a
sequence of eigenvalues $\left\{ \mu _{n}=\left( \frac{n\pi }{a}\right)
^{2}\right\} $ and eigenfunctions 
\begin{equation}
u_{n}=\left\{ \cos \nu _{n}x+\sigma \left( \nu _{n}\right) \sin \nu
_{n}x\right\} \text{, for }n=1,2,...,\infty \text{,}  \tag{1.9}
\end{equation}%
where 
\[
\nu _{n}=\frac{n\pi }{a}\text{, }\sigma \left( \nu _{n}\right) =-\frac{%
\alpha _{2}+\beta _{2}\cos \nu _{n}a}{\alpha _{1}+\beta _{1}\cos \nu _{n}a}=-%
\frac{\alpha _{2}+\beta _{2}\left( -1\right) ^{n}}{\alpha _{1}+\beta
_{1}\left( -1\right) ^{n}}. 
\]%
This family of eigenfunctions is an orthogonal basis in $L^{2}\left(
0,a\right) $ if 
\[
\sigma \left( \nu _{n}\right) \left( \frac{1}{\nu _{n}+\nu _{k}}-\frac{1}{%
\nu _{n}-\nu _{k}}\right) +\sigma \left( \nu _{k}\right) \left[ \frac{1}{\nu
_{n}+\nu _{k}}-\frac{1}{\nu _{k}-\nu _{n}}\right] = 
\]

\[
\frac{2a\left[ k\sigma \left( \nu _{k}\right) -j\sigma \left( \nu
_{n}\right) \right] }{\left( n^{2}-k^{2}\right) \pi }=0\text{ for }n\neq j%
\text{, i.e., }\alpha _{1}+\beta _{1}\left( -1\right) ^{n}\neq 0\text{,} 
\]%
\begin{equation}
\text{ }n\left[ \frac{\alpha _{2}+\beta _{2}\left( -1\right) ^{n}}{\alpha
_{1}+\beta _{1}\left( -1\right) ^{n}}\right] -k\left[ \frac{\alpha
_{2}+\beta _{2}\left( -1\right) ^{k}}{\alpha _{1}+\beta _{1}\left( -1\right)
^{k}}\right] =0\text{ for }n\neq k.  \tag{1.10}
\end{equation}%
For $S=-\partial _{x}^{2}I$ in $L^{2}\left( 0,a;H\right) $ with boundary
conditions $\left( 1.8\right) $, we set $W_{r}=D\left( S^{\frac{r}{2}%
}\right) $. So, $W_{0}=L^{2}\left( 0,a;H\right) $, $W_{1}=W_{a}^{1}\left(
0,a;H\right) $, and 
\[
W_{2}=W^{2}\left( 0,a;H\right) \cap W_{a}^{1}\left( 0,a;H\right) , 
\]%
here 
\[
W_{a}^{1}\left( 0,a;H\right) =\left\{ {}\right. u\text{: }u\in W^{1}\left(
0,a;H\right) \text{, } 
\]%
\[
\alpha _{1}u\left( 0\right) +\beta _{1}u\left( a\right) =0,\text{ }\alpha
_{2}u\left( 0\right) +\beta _{2}u\left( a\right) =0\left. {}\right\} . 
\]

By reasoning as in $\left[ 13\right] $ and method of Hilbert spaces we have
the following generalizasion of B. S. Pavlov theorem $\left[ \text{13}\right]
$:

\textbf{Theorem A}$_{1}$. Let $\Lambda =\left\{ \lambda _{k}\text{: }k\in 
\mathbb{Z}\right\} $ be a countable set in the complex plane and $\left\{
\varphi _{j}\right\} _{j=1}^{\infty }$ is a Riesz basis in a Hilbert space $%
H.$ The family $\left\{ \varphi _{j}\exp \left\{ i\lambda _{k}t\right\}
\right\} $ forms a Riesz basis in \ $L^{2}\left( 0,T;H\right) $ if the
following conditions are satisfied:

(i) $\Lambda $ lies in a strip parallel to the real axis, $\sup\limits_{k\in 
\mathbb{Z}}\left\vert \func{Im}\lambda _{k}\right\vert <\infty $ and is
uniformly discrete (or separated), i.e.%
\begin{equation}
\delta \left( \Lambda \right) \text{:}=\inf\limits_{k\neq n}\left\vert
\lambda _{k}-\lambda _{n}\right\vert >0;  \tag{1.6}
\end{equation}

(ii) there exists an entire function $F$ of exponential type with indicator
diagram of width $T$ and zero set (the generating function of the family fei
ntg on the interval $(0;T)$ such that, for some real $h$, the function $%
\left\Vert F\left( x+ih\right) \right\Vert _{H}^{2}$\ satisfies the
Helson-Szego condition: functions $u$, $\upsilon \in L^{\infty }\left( 
\mathbb{R}\right) $, $\left\Vert \upsilon \right\Vert _{L^{\infty }\left( 
\mathbb{R}\right) }<\frac{\pi }{2}$ can be found such that

\begin{equation}
\left\Vert F\left( x+ih\right) \right\Vert _{H}^{2}=\exp \left\{ u\left(
x\right) +\tilde{\upsilon}\left( x\right) \right\} ;  \tag{1.7}
\end{equation}
here, the map $\upsilon \rightarrow \tilde{\upsilon}$ denotes the Hilbert
transform for bounded functions:%
\[
\tilde{\upsilon}=H\upsilon =\frac{1}{\pi }p.v\dint\limits_{-\infty }^{\infty
}\left[ \frac{1}{x-t}+\frac{t}{t^{2}+1}\right] \upsilon \left( t\right) dt. 
\]

\textbf{Condition 1.1. }Let $A$ be a symmetric operator in a Hilbert space $%
H $ with discrete specrum $\left\{ \lambda _{k}\right\} _{k=1}^{\infty }$
and corresponding eigenfunctions $\left\{ \varphi _{k}\right\}
_{k=1}^{\infty }$. Moreover, assume that $\left\{ \varphi _{k}\right\}
_{k=1}^{\infty }$ is a ortonormal system in $H$.

Our main result is the following:

\textbf{Theorem 1.1.} \ Let $\ $the Condition 1.1. holds and $A$ have the
distinct eigenvalues $\left\{ \lambda _{k}\right\} $, $k\in \mathbb{N}$.
Moreover, assume the following assumptions are satisfied:

(1) $\ $a linear operator $A$ and $b\in H$ such that the system $\left\{
A^{k}b\right\} $, $k\in \mathbb{N}$ is linearly independent in Hilbert space 
$H;$

(2) $\mu _{k}-\mu _{l}\neq \lambda _{i}-\lambda _{j}$ for each $k,l\in 
\mathbb{N},$ $i,j\in \left\{ 1,2,...\right\} $ with $k\neq l$ and $i\neq j;$

(3) $\eta =\alpha _{1}\beta _{2}-\alpha _{2}\beta _{1}\neq 0$ and $\left(
1.5\right) $ is satisfied.

Then the problem $\left( 1.1\right) $ is exactly controllable in $\mathbb{H=}%
W_{0}\times W_{-1}.$

\begin{center}
\bigskip \textbf{2. Proof of Theorem 1.1.}
\end{center}

\textbf{The existence of solutions}. In this section, we use the Fourier
method and apply it to the case, where the operator $A$ has distinct
eigenvalues. On the assumptions of Theorem 1.1 we denote $\left\{ \lambda
_{k}\right\} $ to be the family of eigenvectors of $A$ with corresponding
eigenfunctions $\left\{ \varphi _{k}\right\} _{k=1}^{\infty }$. We denote by 
$\ \left( .,.\right) =\left( .,.\right) _{H}$, $\langle .,.\rangle $ the
inner product in Hilbert spaces $H$ and $L^{2}\left( \Omega ;H\right) $,
respectively. So conjucate operator $A^{\ast }$ has eigenvalues $\left\{ 
\bar{\lambda}_{k}\right\} $ and eigenvectors $\left\{ \psi _{k}\right\}
_{k=1}^{\infty }$ with 
\[
\left( \varphi _{i},\psi _{j}\right) _{H}=\delta _{ij}. 
\]

Let we give some lemmas for proving of Theorem 1.1.

\bigskip \textbf{Lemma 2.1. }Let the Condition 1.1. holds. Suppose that $A$
have distinct eigenvalues $\lambda _{1}$, $\lambda _{2}$, . . .$\lambda
_{n},...$. Assuming that the assumption (1) is satisfied. Then eigenvectors $%
\left\{ \varphi _{k}\right\} _{k=1}^{\infty }$ and $\left\{ \psi
_{k}\right\} _{k=1}^{\infty }$ may be chosen such that $\langle b,\psi
_{j}\rangle =1$\ while\ $\langle \varphi _{i},\psi _{j}\rangle =\delta
_{ij}. $

\textbf{Proof.} \ We first claim that $\left( b,\psi _{j}\right) \neq 0$.
Indeed, if there exists $i\in \mathbb{N}$ such that $\left( b,\psi
_{i}\right) =0$, then for all $m\in \mathbb{N},$%
\[
\left( A^{m}b,\psi _{i}\right) =\left( b,\left( A^{\ast }\right) ^{m}\psi
_{i}\right) =\left( b,\left( \bar{\lambda}_{i}\right) ^{m}\psi _{i}\right)
=\lambda _{i}^{m}\left( b,\psi _{i}\right) =0. 
\]

This implies that the system $\left\{ A^{k}b\right\} $, $k\in \mathbb{N}$ is
linearly dependent in Hilbert space $H$, which is a contradiction the
ussumption (1). Hence, we can construct the sets $\left\{ \tilde{\varphi}%
_{i}\right\} $, $\left\{ \tilde{\psi}_{i}\right\} $, $i\in \mathbb{N}$, where

\[
\tilde{\varphi}_{i}=\left( b,\psi _{i}\right) \varphi _{i}\text{, }\tilde{%
\psi}_{i}=\frac{\psi _{i}}{\left( b,\psi _{i}\right) }. 
\]

It then follows that $\left( b,\tilde{\psi}_{i}\right) =1$ for $i\in \mathbb{%
N}$ and $\left( \tilde{\varphi}_{i},\tilde{\psi}_{j}\right) =\delta _{ij}$
for $i$, $j\in \mathbb{N}$. \ So we may assume that $\left( b,\psi
_{i}\right) =1.$

\bigskip Let us define $\Phi _{nk}\left( x\right) =u_{n}\left( x\right)
\varphi _{k}$, where $\left\{ u_{n}\right\} $ is a system defined by $\left(
1.9\right) $. Then $\left\{ \Phi _{nk}\left( x\right) \right\} $, $n,k\in 
\mathbb{N}$ is a basis in $L^{2}\left( 0,a;H\right) $ with biorthogonal
family $\left\{ \Psi _{nk}\left( x\right) \right\} =\left\{ \bar{u}%
_{n}\left( x\right) \psi _{k}\right\} ,$ where

\[
\dint\limits_{0}^{a}u_{k}\left( x\right) \bar{u}_{n}\left( x\right)
dx=\delta _{kn}, 
\]%
i.e, for example, 
\begin{equation}
\bar{u}_{n}\left( x\right) =\gamma _{n}u_{n}\left( x\right) \text{, }\gamma
_{n}=\frac{a}{2}\left( 1+\sigma ^{2}\nu _{n}\right) .  \tag{2.1}
\end{equation}

We then can represent the solution $u$ of the problem $\left( 1.1\right) $
in the form of the series%
\begin{equation}
u\left( x,t\right) =\dsum\limits_{n,k}a_{nk}\left( t\right) \Phi _{nk}\left(
x\right)  \tag{2.2}
\end{equation}%
and set 
\begin{equation}
\upsilon \left( x,t\right) =g\left( t\right) \Psi _{kl}\left( x\right) 
\tag{2.3}
\end{equation}%
for some $j,l\in \mathbb{N}$, where $g\left( t\right) \in C_{0}^{2}\left(
0,T\right) $ such that 
\begin{equation}
\dint\limits_{0}^{T}\left[ \left( \left( -1\right) ^{k}u_{x}\left(
a,t\right) -u_{x}\left( 0,t\right) ,\psi _{l}\right) _{H}\right] g\left(
t\right) =0\text{ for }k\text{, }l\in \mathbb{N}.  \tag{2.4}
\end{equation}

Let the function expressed $\left( 2.2\right) $ is a solution of $\left(
1.1\right) $. Then from $\left( 1.1\right) $ we get 
\[
\dint\limits_{0}^{T}\dint\limits_{0}^{a}\left( u_{tt}-u_{xx}+Au,\upsilon
\right) _{H}dxdt=\dint\limits_{0}^{T}\dint\limits_{0}^{a}\left( u,\upsilon
_{tt}-\upsilon _{xx}+A^{\ast }\upsilon \right) _{H}dxdt+ 
\]%
\begin{equation}
\dint\limits_{0}^{a}\left[ \left( u_{t},\upsilon \right) _{H}-\left(
u,\upsilon _{t}\right) _{H}\right] _{t=0}^{T}dx-\dint\limits_{0}^{T}\left[
\left( u_{x},\upsilon \right) _{H}-\left( u,\upsilon _{x}\right) _{H}\right]
_{x=0}^{a}dt.  \tag{2.5}
\end{equation}

By $\left( 1.9\right) $, $\left( 2.1\right) $ and $\left( 2.2\right) $ we
have%
\[
\bar{u}_{k}\left( 0\right) =\gamma _{k}\text{, }\bar{u}_{k}^{\left( 1\right)
}\left( 0\right) =\nu _{k}\gamma _{k}\sigma \left( \nu _{k}\right) , 
\]%
\begin{equation}
\bar{u}_{k}\left( a\right) =\gamma _{k}\left( -1\right) ^{k}\text{, }\bar{u}%
_{k}^{\left( 1\right) }\left( a\right) =\left( -1\right) ^{k}\nu _{k}\gamma
_{k}\sigma \left( \nu _{k}\right) \text{. }  \tag{2.6}
\end{equation}

\bigskip Since $\Psi _{nk}\left( x\right) =\bar{u}_{n}\left( x\right) \psi
_{k}$, from $\left( 2.3\right) $-$\left( 2.6\right) $ by taking the nonlocal
problem $\left( 1.1\right) $, by using Lemma 2.1 and in view of $g\in
C_{0}^{2}\left( 0,T\right) $ we obtain

\[
\dint\limits_{0}^{T}\dint\limits_{0}^{a}\left( u_{tt}-u_{xx}+Au,\upsilon
\right) _{H}dxdt=\dint\limits_{0}^{T}\dint\limits_{0}^{a}\left( u,\upsilon
_{tt}-\upsilon _{xx}+A^{\ast }\upsilon \right) dxdt= 
\]%
\[
\dint\limits_{0}^{T}\dint\limits_{0}^{a}\left( u,\Psi _{kl}g^{\left(
2\right) }-g\Psi _{kl}\partial _{x}^{2}u_{k}+\bar{\lambda}_{k}g\Psi
_{kl}\right) dxdt+ 
\]%
\[
\dint\limits_{0}^{T}\left[ \left( u_{x}\left( a,t\right) ,\Psi _{kl}\left(
a\right) \right) +\frac{\alpha _{2}}{\eta }f\left( t\right) \left( b,\Psi
_{kl}^{\prime }\left( a\right) \right) \right] g\left( t\right) dt- 
\]%
\[
\dint\limits_{0}^{T}\left[ \left( u_{x}\left( 0,t\right) ,\Psi _{kl}\left(
0\right) \right) -\frac{\beta _{2}}{\eta }f\left( t\right) \left( b,\Psi
_{kl}^{\prime }\left( 0\right) \right) \right] g\left( t\right) dt= 
\]%
\[
\dint\limits_{0}^{T}\dint\limits_{0}^{a}\left( u,\Psi _{kl}g^{\left(
2\right) }+\nu _{k}^{2}g\Psi _{kl}+\bar{\lambda}_{l}\Psi _{kl}\right)
g\left( t\right) dxdt+ 
\]

\[
\frac{\gamma _{k}\nu _{k}\sigma \left( \nu _{k}\right) }{\eta }\left( \beta
_{2}-\alpha _{2}\right) \dint\limits_{0}^{T}\left( b,\psi _{l}\right)
f\left( t\right) g\left( t\right) dt= 
\]%
\[
\dint\limits_{0}^{T}a_{kl}\left[ g^{\left( 2\right) }+\left( \nu _{k}^{2}+%
\bar{\lambda}_{l}\right) \right] g\left( t\right) dt+\varkappa
_{k}\dint\limits_{0}^{T}f\left( t\right) g\left( t\right) dt= 
\]%
\[
\dint\limits_{0}^{T}\left[ a_{kl}^{\left( 2\right) }+\left( \nu _{k}^{2}+%
\bar{\lambda}_{l}\right) a_{kl}\right] g\left( t\right) dt-\varkappa
_{k}\dint\limits_{0}^{T}f\left( t\right) g\left( t\right) dt=0, 
\]%
where by Lemma 1, $\left( b,\psi _{l}\right) _{H}=1$ and 
\[
\varkappa _{k}=\frac{\gamma _{k}\nu _{k}\sigma \left( \nu _{k}\right) }{\eta 
}\left( \alpha _{2}-\beta _{2}\right) . 
\]%
Thus we obtain the equations%
\begin{equation}
a_{kl}^{\prime \prime }+\left( \nu _{k}^{2}+\bar{\lambda}_{l}\right)
a_{kl}^{\prime }=\varkappa _{k}f\left( t\right)  \tag{2.7.}
\end{equation}%
with initial conditions 
\begin{equation}
a_{kl}\left( 0\right) =a_{kl}^{\prime }\left( 0\right) =0.  \tag{2.8}
\end{equation}

We assume 
\begin{equation}
\omega _{kl}=\left( \nu _{k}^{2}+\bar{\lambda}_{l}\right) ^{\frac{1}{2}}\neq
0\text{, }\nu _{k}=\frac{k\pi }{a}.  \tag{2.9}
\end{equation}

\bigskip We can set the following

\textbf{Proposition 2.1}. Let $k\in \mathbb{K=}\left\{ \pm 1\text{, }\pm
2,...\right\} $ and $1\leq m,l\leq n$ with $m\neq l$. Provided the condition
(2) of Theorem 1.1, we have the following:

(1) $\left\vert \omega _{kl}\right\vert +1\asymp k;$

(2) $\left\vert \omega _{kl}-\omega _{km}\right\vert \asymp k^{-1};$

(3) For $k$ fixed, the points $\omega _{kl}$ are asymptotically close, i.e.,
these points lie inside an interval whose length tends to zero as $k$ tends
to infinity.

Let 
\[
X=L^{2}\left( 0,a;H\right) \text{, }W^{s}=W^{s,2}\left( 0,a;H\right) \text{, 
}Y=X\times W^{s}. 
\]%
By following $\left[ \text{9, Theorem 2.1}\right] $ we have:

\textbf{Theorem 2.1. }Let $\ $the Condition 1.1. holds and $A$ have distinct
eigenvalues $\lambda _{1}$, $\lambda _{2}$, . . .$\lambda _{n},...$. Then\
for any $f\in L^{2}\left( 0,T\right) $ there exists a unique generalized
solution $u=u^{f}\left( x,t\right) $ of the problem $\left( 1.1\right) $
such that $\left( u^{f},u_{t}^{f}\right) \in C\left( \left[ 0,T\right]
;Y\right) $ and 
\[
\left\Vert \left( u^{f},u_{t}^{f}\right) \right\Vert _{C\left( \left[ 0,T%
\right] ;Y\right) }\prec \left\Vert f\right\Vert _{L^{2}\left( 0,T\right) }. 
\]

\textbf{Proof. }The solution of $\left( 2.7\right) -\left( 2.8\right) $ is
geven by the formula 
\begin{equation}
a_{kl}\left( t\right) =\varkappa _{k}\dint\limits_{0}^{t}f\left( s\right) 
\frac{\sin \omega _{kl}\left( t-s\right) }{\omega _{kl}}ds.  \tag{2.10}
\end{equation}

By differentiating we have%
\begin{equation}
a_{kl}^{\prime }\left( t\right) =\varkappa _{k}\dint\limits_{0}^{t}f\left(
s\right) \cos \omega _{kl}\left( t-s\right) ds.  \tag{2.11}
\end{equation}

We now introduce the coefficients 
\begin{equation}
c_{kl}\left( t\right) =i\omega _{kl}a_{kl}\left( t\right) +a_{kl}^{\prime
}\left( t\right) .  \tag{2.12}
\end{equation}

\begin{center}
\bigskip
\end{center}

\bigskip Now, we define 
\[
\omega _{-kl}=-\omega _{kl}\text{, }a_{-kl}=a_{kl}\text{, and }%
a_{-kl}^{\prime }=a_{kl}^{\prime }\text{ for }k\in \mathbb{K}\text{, }l\in 
\mathbb{N} 
\]%
and rewrite $\left( 2.10\right) $, $\left( 2.11\right) $ in the exponential
form, we get 
\begin{equation}
C_{kl}\left( t\right) =\varkappa _{k}\dint\limits_{0}^{t}f\left( s\right)
\exp i\left( \omega _{kl}\left( t-s\right) \right) ds.  \tag{2.13}
\end{equation}

Taking into account that $\Phi _{nj}$ forms a Riesz basis in \ $L^{2}\left(
0,a;H\right) $ and Proposition 1 property (1), by $\left[ {}\right] $\ we
conclude that%
\begin{equation}
\dsum\limits_{k\in \mathbb{K}}\frac{\left[ C_{kl}\left( t\right) \right] ^{2}%
}{k^{2}}\asymp \left\Vert u\left( .,t\right) \right\Vert _{L^{2}\left(
0,a;H\right) }^{2}+\left\Vert u_{t}\left( .,t\right) \right\Vert
_{H^{-1}\left( 0,a;H\right) }^{2}.  \tag{2.14}
\end{equation}

On the other hand, from the explicit form for !kl, it follows that for any $%
T>0$, the family $\left\{ \exp i\left( \omega _{kl}t\right) \right\} $ is
either a finite union of Riesz sequences if $T<2na$ or a Riesz sequence in $%
L^{2}\left( 0,T\right) $ if $T>2na$ (see $\left[ 8\right] $ Section II.4).
We recall that a Riesz sequence is a Riesz basis in the closure of its
linear span. Therefore, from $\left( 2.13\right) $ it follows that for every
fixed $t>0,$%
\begin{equation}
\dsum\limits_{k,l}\frac{\left\vert C_{kl}\left( t\right) \right\vert ^{2}}{%
k^{2}}\prec \left\Vert f\right\Vert _{L^{2}\left( 0,T\right) }^{2}. 
\tag{2.15}
\end{equation}

It can be shown that the series in $\left( 2.15\right) $ is uniformly
convergent by the Weierstrass criterion for uniform convergence. And by the
uniform limit theorem, we obtain%
\[
\dsum\limits_{k,l}\frac{\left\vert C_{kl}\left( t+h\right) -C_{kl}\left(
t\right) \right\vert ^{2}}{k^{2}}\rightarrow 0\text{ when }h\rightarrow 0. 
\]

\textbf{2.2. Controllability results\bigskip . }In this section we will
prove Theorem 1.1. Let 
\begin{equation}
\alpha _{kl}=C_{kl}\left( T\right) \left( \frac{2k}{\pi }\exp \left\{
i\omega _{kl}T\right\} \right) ^{-1}  \tag{2.16}
\end{equation}

and rewrite $\left( 2.13\right) $ for $t=T$ in the form 
\begin{equation}
\alpha _{kl}=\left( f,e_{kl}\right) _{L^{2}\left( 0,a\right) },  \tag{2.17}
\end{equation}%
where 
\[
e_{kl}=e_{kl}\left( t\right) =\exp \left\{ i\omega _{kl}t\right\} . 
\]

\bigskip We note that 
\[
\dsum\limits_{k,l}\left\vert \alpha _{kl}\right\vert ^{2}\asymp
\dsum\limits_{k,l}\frac{\left\vert C_{kl}\left( T\right) \right\vert ^{2}}{%
k^{2}}. 
\]

\bigskip For any $T>0$, the family $\left\{ e_{kl}\right\} $ is not a Riesz
basis as a result of Proposition 1 property (3). Therefore, we need to use
the so-called exponential divided di erences (EDD). EDD were introduced in $%
[5]$ and $[6]$ for families of exponentials whose exponents are close, that
is, the di erence between exponents tends to zero. Under precise
assumptions, the family of EDD forms a Riesz sequence in $L^{2}\left(
0,T\right) $. For each fixed $k$, we define 
\[
\tilde{e}_{kl}\text{:}=\left[ \omega _{k1}\right] =\exp \left\{ i\omega
_{k1}t\right\} , 
\]%
and for $2\leq l\leq d$, we define 
\[
\tilde{e}_{kl}\text{:}=\left[ \omega _{k1},\omega _{k2},...,\omega _{kl}%
\right] =\dsum\limits_{j=1}^{l}\frac{\exp \left\{ i\omega _{kj}t\right\} }{%
\dprod\limits_{n\neq j}\left( \omega _{kj}-\omega _{nj}\right) }. 
\]

Under Condition (ii) of our theorem, we are able to use this formula for
divided differences in place of the formula for generalized divided
differences (see e.g. $[10]$). From asymptotics theory and the explicit
formula for $\omega _{kl}$, it follows that the generating function of the
family of EDD $\left\{ \tilde{e}_{kl}\right\} $ is a sine-type function (see 
$[9,10]$). Hence, the family of EDD $\left\{ \tilde{e}_{kl}\right\} $ forms
a Riesz sequence in $L^{2}\left( 0,T\right) $. We then define%
\[
\tilde{\alpha}_{kl}=<f,\tilde{e}_{kl}>. 
\]

\bigskip Since $\left\{ \tilde{e}_{kl}\right\} $ is a Riesz sequence, $%
\left\{ \tilde{\alpha}_{kl}\text{:}f\in L^{2}\left( 0,T\right) \right\}
=l_{2}$, i.e. any sequence from $l_{2}$ can be obtained by a function $f\in
L^{2}\left( 0,T\right) $ and the family $\left\{ \tilde{e}_{kl}\right\} $.
Proposition 1 property (2) implies that $\left\vert \tilde{\alpha}%
_{kl}\right\vert \prec k^{d-1}\left\vert \alpha _{kl}\right\vert $. Then by
reasoning as in $\left[ \text{7, \S\ 2.2}\right] $ we obtain the assertion
of Theorem 1.1.

\begin{center}
\bigskip \textbf{3}. \textbf{Application}
\end{center}

\bigskip \textbf{3.1. Boundary controllability for infinite many system of
wave equations. }Consider the problem 1.3. Let%
\[
\text{ }l_{2}\left( N\right) =\left\{ \text{ }u=\left\{ u_{j}\right\} ,\text{
}j=1,2,...N,\left\Vert u\right\Vert _{l_{2}\left( N\right) }=\left(
\sum\limits_{j=1}^{N}\left\vert u_{j}\right\vert ^{2}\right) ^{\frac{1}{2}%
}<\infty \right\} , 
\]%
where $N\in \mathbb{N}$ (see $\left[ \text{15, \S\ 1.18}\right] .$ Let $A$
be the operator in $l_{2}\left( N\right) $ defined by%
\begin{equation}
\text{ }A=\left[ a_{jm}\right] ,\text{ }a_{jm}=b_{j}2^{\sigma m},\text{ }%
m,j=1,2,...,N,\text{ }D\left( A\right) =\text{ }l_{2}^{\sigma }\left(
N\right) =  \tag{3.1}
\end{equation}

\[
\left\{ \text{ }u=\left\{ u_{j}\right\} ,\text{ }j=1,2,...N,\left\Vert
u\right\Vert _{l_{2}^{\sigma }\left( N\right) }=\left(
\sum\limits_{j=1}^{N}2^{\sigma j}u_{j}^{2}\right) ^{\frac{1}{2}}<\infty
\right\} ,\text{ }\sigma >0. 
\]

Let $l_{2}=l_{2}\left( N\right) $, $b=\left\{ b_{m}\right\} $, $m=1,2,...,N$
and \ 
\[
X=L^{2}\left( 0,a;l_{2}\right) \text{, }X^{s}=H^{s,2}\left( 0,a;l_{2}\right)
. 
\]

\ From Theorem 1.1 we obtain:

\textbf{Theorem 3.1.} \ Suppose that:

(1) $a_{jm}\in \mathbb{R}$, $a_{jm}=a_{mj}\ $and $A$ have the distinct
eigenvalues $\left\{ \lambda _{k}\right\} $, $k\in \mathbb{N}$;

(2) $A$ and $b\in H$ such that the system $\left\{ A^{k}b\right\} $, $k\in 
\mathbb{N}$ is linearly independent in $l_{2};$

(3) $\mu _{k}-\mu _{l}\neq \lambda _{i}-\lambda _{j}$ for each $k,l\in 
\mathbb{N},$ $i,j\in \left\{ 1,2,...\right\} $ with $k\neq l$ and $i\neq j;$

(4) $\alpha _{1}\beta _{2}-\alpha _{2}\beta _{1}\neq 0$ and $\left(
1.10\right) $ is satisfied.

Then the problem $\left( 1.3\right) $ is exactly controllable in $\mathbb{H=}%
X\times X^{-1}.$

\textbf{Proof. }\ It is clear to see that the operator $A$ defined in
Hilbert space $l_{2}$ defined by $\left( 3.1\right) $ is symmetric. By
assumptions (1)-(4) all conditions of Theorem 1.1. is satisfied, i.e. we
obtain the assertion.

\bigskip \textbf{3.2. Boundary controllability for degenerate wave
equations. }Consider the problem 1.4-1.7. Let

Let \ $b\in L^{2}\left( 0,1\right) $ and 
\[
Y=L^{2}\left( 0,a;L^{2}\left( 0,1\right) \right) =L^{2}\left( \left(
0,a\right) \times \left( 0,1\right) \right) \text{, }Y^{s}=H^{s,2}\left(
0,a;L^{2}\left( 0,1\right) \right) . 
\]

Consider the operator in $L^{2}\left( 0,1\right) $ defined by 
\begin{equation}
D\left( A\right) =W_{\gamma }^{\left[ 2\right] ,2}\left( 0,1\right) \text{, }%
Au=\left( a_{1}\frac{d^{\left[ 2\right] }u}{dy^{2}}+a_{2}\frac{d^{\left[ 1%
\right] }u}{dy}\right) .  \tag{3.2}
\end{equation}

From Theorem 1.1 we obtain:

\textbf{Theorem 3.2.} \ Assume that:

(1) $a_{1}$ is positive continious and $a_{2}$ a bounded functions and $A$
is a symmetric operator in $L^{2}\left( 0,1\right) $\ having the distinct
eigenvalues;

(2) $A$ and $B\in H$ such that the system $\left\{ A^{k}b\right\} $, $k\in 
\mathbb{N}$ is linearly independent in $l_{2}$ and $0\leq \gamma <\frac{1}{2}%
;$

(3) $\mu _{k}-\mu _{l}\neq \lambda _{i}-\lambda _{j}$ for each $k,l\in 
\mathbb{N},$ $i,j\in \left\{ 1,2,...\right\} $ with $k\neq l$ and $i\neq j;$

(4) $\alpha _{1}\beta _{2}-\alpha _{2}\beta _{1}\neq 0$, $\nu _{1}\delta
_{2}-\nu _{2}\delta _{1}$ and $\left( 1.10\right) $ is satisfied.

Then the problem $\left( 1.3\right) $ is exactly controllable in $\mathbb{H=}%
X\times Y^{-1}.$

\textbf{Proof. }\ By $\left[ \text{14, Theorem 3.3}\right] $, the operator $%
A $ defined in $L^{2}\left( 0,1\right) $ by $\left( 3.2\right) $ have a
discrete specrum. By assumptions (1)-(4) all conditions of Theorem 1.1. is
satisfied, i.e. we obtain the assertion.

\begin{center}
\bigskip

\textbf{References}
\end{center}

\bigskip

[1] F. Alabau-Boussouira, A two-level energy method for indrect boundary
observability and controllability of weakly coupled hyperbolic systems, SIAM
J. Control Optim., 42 (2003), 871-906.

[2] F. Alabau-Boussouira, Insensitizing exact controls for the scalar wave
equation and exact controllability of 2-coupled cascade systems of PDE's by
a single control, Math. Control Signals Systems, 26 (2014), 1-46.

[3] F. Alabau-Boussouira and M. L eautaud, Indirect controllability of
locally coupled systems under geometric conditions, C. R. Acad. Sci. Paris,
349 (2011), 395-400.

[4] F. Ammar-Kohdja, A. Benabdallah, M. Gonz alez-Burgos and L. de Teresa,
The Kalman condition for the boundary controllability of coupled parabolic
systems. bounds on biorthogonal families to complex matrix exponentials,
JMPA, 96 (2011), 555\{590, https://doi.org/10. 1016/j.matpur.2011.06.005.

$\left[ 5\right] $ \ S. A. Avdonin and S. A. Ivanov, Exponential Riesz bases
of subspaces and divided di erences, St. Petersburg Mathematical Journal, 13
(2002), 339\{351.

[6] S. Avdonin and W. Moran, Ingham type inequalities and Riesz bases of
subspaces and divided di erences, Int. J. Appl. Math. Compt. Sci., 11
(2001), 803-820.

[7] S. Avdonin and L. de Teresa, The Kalman Condition for the Boundary
Controllability of Coupled 1-d Wave Equations, Evol. Equat. and Cont.Theory,
1(9) (2020), 255-273.

[8] S. A. Avdonin and S. A. Ivanov, Families of Exponentials: The Method of
Moments in Con- trollability Problems for Distributed Parameter Systems,
Cambring University Press, 1995.

[9] S. A. Avdonin, J. Park 1, L. de Teresa, The Kalman condition for the
boundary controllability of coupled 1-d wave equations, Evolution equations
and control theory, (1) 9 (2020), 255-273.

[10] S. Avdonin, A. Choque and L. de Teresa, Exact boundary controllability
results for two coupled 1-d hyperbolic equations, Int. J. Appl. Math.
Comput. Sci., 23 (2013), 701-710, https://doi.org/10.2478/amcs-2013-0052.

[11] A. Bennour, F. Ammaar Khodja and D. Tenious, Exact and approximate
controllability of coupled one-dimensional hyperbolic equations, Ev. Eq. and
Cont. Teho., 6 (2017), 487\{516.

$\left[ 12\right] $ R. E. Kalman, P. L. Palb and M. A. Arbib, Topics in
Mathematical Control Theory, New York-Toronto, Ont.-London, 1969.

$\left[ 13\right] $ B. S. Pavlov, Basicity of exponential system and
Muckenhoupt condition, Doklady Akad. Nauk. SSSR. 247(1) (1979), 37-40
(Russian); English transl. in Soviet Math. Dokl. 20: 655-9.

$\left[ 14\right] $ V. B. Shakhmurov, Linear and nonlinear abstract
differential equations with small parameters, Banach J. Math. Anal. 10(1)
(2016), 147--168.

$\left[ 15\right] $ Interpolation theory, Function spaces, Differential
operators, North-Holland, Amsterdam, 1978.

\end{document}